\documentclass[12pt,leqno,draft]{article}
\usepackage{amsfonts}
\usepackage{amsmath, amsthm, amsfonts, amssymb, color}
\usepackage{mathrsfs}
\usepackage{color}
\theoremstyle{definition}

\newcommand{\scr}[1]{\mathscr #1}
\definecolor{wco}{rgb}{0.5,0.2,0.3}

\numberwithin{equation}{section} \theoremstyle{remark}

\newcommand{\ua}{\uparrow}

\title{{\bf
 Distribution Dependent Birth-Death Processes: $\W_p$-Estimate, Ergodicity and Propagation of Chaos}\footnote{Supported in part by the National Key R\&D Program of China (2022YFA1006000), 
 NSFC(12531007) and the State Key Laboratory of Synthetic Biology in Tianjin University.} }
\author{
{\bf     Feng-Yu Wang, Yi Zhao }\\
\footnotesize{  Center for Applied Mathematics and KL-AAGDM, Tianjin
University, Tianjin 300072, China} \\
\footnotesize{   wangfy@tju.edu.cn, zhaoyibnu@mail.bnu.edu.cn}\\
}
\begin{document}
\allowdisplaybreaks
\def\R{\mathbb R}  \def\ff{\frac} \def\ss{\sqrt} \def\B{\mathbf
B} \def\W{\mathbb W}
\def\N{\mathbb N} \def\kk{\kappa} \def\m{{\bf m}}
\def\ee{\varepsilon}\def\ddd{D^*}
\def\dd{\delta} \def\DD{\Delta} \def\vv{\varepsilon} \def\rr{\rho}
\def\<{\langle} \def\>{\rangle} \def\GG{\Gamma} \def\gg{\gamma}
  \def\nn{\nabla} \def\pp{\partial} \def\E{\mathbb E}
\def\d{\text{\rm{d}}} \def\bb{\beta} \def\aa{\alpha} \def\D{\scr D}
  \def\si{\sigma} \def\ess{\text{\rm{ess}}}
\def\beg{\begin} \def\beq{\begin{equation}}  \def\F{\scr F}
\def\Ric{\text{\rm{Ric}}} \def\Hess{\text{\rm{Hess}}}
\def\e{\text{\rm{e}}} \def\ua{\underline a} \def\OO{\Omega}  \def\oo{\omega}
 \def\tt{\tilde} \def\Ric{\text{\rm{Ric}}}
\def\cut{\text{\rm{cut}}} \def\P{\mathbb P} \def\ifn{I_n(f^{\bigotimes n})}
\def\C{\scr C}      \def\aaa{\mathbf{r}}     \def\r{r}
\def\gap{\text{\rm{gap}}} \def\prr{\pi_{{\bf m},\varrho}}  \def\r{\mathbf r}
\def\Z{\mathbb Z} \def\vrr{\varrho} \def\ll{\lambda}
\def\L{\scr L}\def\Tt{\tt} \def\TT{\tt}\def\II{\mathbb I}
\def\i{{\rm in}}\def\Sect{{\rm Sect}}  \def\H{\mathbb H}
\def\M{\scr M}\def\Q{\mathbb Q} \def\texto{\text{o}} \def\LL{\Lambda}
\def\Rank{{\rm Rank}} \def\B{\scr B} \def\i{{\rm i}} \def\HR{\hat{\R}^d}
\def\to{\rightarrow}\def\l{\ell}\def\iint{\int}
\def\EE{\scr E}\def\Cut{{\rm Cut}}
\def\A{\scr A} \def\Lip{{\rm Lip}}
\def\BB{\scr B}\def\Ent{{\rm Ent}}\def\L{\scr L}
\def\R{\mathbb R}  \def\ff{\frac} \def\ss{\sqrt} \def\B{\mathbf
B}
\def\N{\mathbb N} \def\kk{\kappa} \def\m{{\bf m}}
\def\dd{\delta} \def\DD{\Delta} \def\vv{\varepsilon} \def\rr{\rho}
\def\<{\langle} \def\>{\rangle} \def\GG{\Gamma} \def\gg{\gamma}
  \def\nn{\nabla} \def\pp{\partial} \def\E{\mathbb E}
\def\d{\text{\rm{d}}} \def\bb{\beta} \def\aa{\alpha} \def\D{\scr D}
  \def\si{\sigma} \def\ess{\text{\rm{ess}}}
\def\beg{\begin} \def\beq{\begin{equation}}  \def\F{\scr F}
\def\Ric{\text{\rm{Ric}}} \def\Hess{\text{\rm{Hess}}}
\def\e{\text{\rm{e}}} \def\ua{\underline a} \def\OO{\Omega}  \def\oo{\omega}
 \def\tt{\tilde} \def\Ric{\text{\rm{Ric}}}
\def\cut{\text{\rm{cut}}} \def\P{\mathbb P} \def\ifn{I_n(f^{\bigotimes n})}
\def\C{\scr C}      \def\aaa{\mathbf{r}}     \def\r{r}
\def\gap{\text{\rm{gap}}} \def\prr{\pi_{{\bf m},\varrho}}  \def\r{\mathbf r}
\def\Z{\mathbb Z} \def\vrr{\varrho} \def\ll{\lambda}
\def\L{\scr L}\def\Tt{\tt} \def\TT{\tt}\def\II{\mathbb I}
\def\i{{\rm in}}\def\Sect{{\rm Sect}}  \def\H{\mathbb H}
\def\M{\scr M}\def\Q{\mathbb Q} \def\texto{\text{o}} \def\LL{\Lambda}
\def\Rank{{\rm Rank}} \def\B{\scr B} \def\i{{\rm i}} \def\HR{\hat{\R}^d}
\def\to{\rightarrow}\def\l{\ell}\def\BB{\mathbb B}
\def\8{\infty}\def\I{1}\def\U{\scr U} \def\n{{\mathbf n}}\def\v{V}\def\LL{{\bf L}}
\maketitle

\begin{abstract}  For a class of time inhomogenous  distribution dependent birth-death processes, we derive the well-posedness,   $\W_p$-estimate, exponential ergodicity, and uniform in time propagation of chaos.
These extend the corresponding results derived for distribution dependent SDEs and mean field particle systems. As preparation, a criterion on the well-posedness of inhomogenous jump process is presented 
in the end of the paper, which should be interesting by itself.   
 \end{abstract} \noindent

\section{Introduction}

Due to its solid background and wide applications, the distribution dependent SDE (DDSDE for short) initiated by McKean \cite{mckean1966} have been intensively investigated, and  plentiful results have been derived 
for the well-posedness, regularities, ergodicity, propagation of chaos, see the monographs \cite{sni91, CD, WR} and references therein. Since the solution of DDSDE is Markov in time but nonlinear in initial distribution, 
it is called   nonlinear Markov process, see  \cite{2010} for a literature of  nonlinear Markov processes  and applications. 

To describe stochastic systems with discontinuous noise, the theory of Markov jump processes has been developed.  This type processes are known as $Q$-processes or continuous time Markov chains when the state space is discrete, and are called $q$-processes in general, see \cite{Ch92} and references therein.  
  In particular, for three specific Markov chains, namely birth-death processes, single birth processes and single death processes, sharp criteria have been derived for different type  ergodicity properties,  see  \cite{Ch92,chen14,zyh18}.  
  
Comparing with  DDSDEs, the study on distribution dependent Markov chains is rather limited. 
  Existing results focus on   expectation dependent Markov chains, i.e. the transition rate depends on the state and expectation of the process,   first introduced by Nicolis and Prigogine \cite{oth82}   to describe a class of nonlinear pure jump processes arising from chemistry, physics and biology, and    further studied  in     \cite{shiga85,feng1992, fs94,fs94ptrf,nonjap} for  the well-posedness,   stationary distributions and     large deviations.   
However,    to our best knowledge, rare is known  on general distribution dependent Markov chains, and the propagation of chaos is open  even in the expectation dependent case.  
      
In this paper, we   study  the well-posedness, exponential ergodicity and propagation of chaos for distribution dependent birth-death processes, and the argument we develop might be  extended to general distribution dependent  jump processes in the future.

\

Let $\Z_+$ be the set of nonnegative integer numbers, let   $\B(\Z_+)$  ($\B_b(\Z_+)$) be the space of all real (bounded) functions  on $\Z_+$, and let  $\scr P$ be the space of probability measures on $\Z_+$ equipped with the weak topology, which is equivalent to the strong topology as all measurable functions on $\Z_+$ are continuous. For a random variable $\xi$ on $\Z_+$, let $\L_\xi\in \scr P$ denote  its  distribution.

 For any $p \in [1,\infty)$,  let
$$\scr P_p:= \Big\{\mu\in \scr P:\ \|\mu\|_p := \mu(|\cdot|^p)^{\ff 1 p}<\infty\Big\}.$$
It  is a Polish space under the $p$-Wasserstein distance 
$$ \W_p(\mu,\nu):= \inf_{\pi\in\mathscr{C}(\mu,\nu)}\left(\int_{\mathbb{Z}_{+}\times\mathbb{Z}_{+}}|x-y|^{p}\pi(\d x,\d y)\right)^{\frac{1}{p}}, $$
where $\mathscr{C}(\mu,\nu)$ is  the set of all couplings of $\mu$ and $\nu$.

Consider the following measurable maps
$$a, b:  [0,\infty)\times \Z_+\times\scr P_1\to [0,\infty),\ \ a_t(0,\mu)=0.$$
   Let 
$$\LL^{a,b}:= \big\{L_{t,\mu}^{a,b}: t\ge 0,\mu\in \scr P_1\big\}$$ be a family of linear operators on $\B(\Z_+)$ defined by
$$L_{t,\mu}^{a,b} f(i):= a_t(i,\mu)\big[f(i-1)-f(i)\big]+ b_t(i,\mu)\big[f(i+1)-f(i)\big],\ \  i\in \Z_+,\  f\in\B(\Z_+),$$ where
$b_t(i,\mu)$ and $a_t(i,\mu)$ stand for the birth and death rates  of the system at time $t$  with state $i$ and  distribution $\mu$. 
The  distribution dependent birth-death process generated by $\LL^{a,b}$ is defined as follows.  

 \beg{defn}\label{def} Let $\gg \in \scr P$. We call $\P^\gg$  
  a  distribution dependent  birth-death process generated by $\LL^{a,b}$, denoted by  $\LL^{a,b}$-process for simplicity,     if   it is a probability measure on the space $D$ consisting of c$\grave{a}$dl$\grave{a}$g
  paths from $[0,\infty) $ to $\Z_+$, such that the time-marginal distributions 
  $$P_t^*\gg(\cdot):=\P^\gg(\{\oo\in D:\ \oo_t\in \cdot\}),\ \ \ t\ge 0$$ satisfy the following conditions: 
  \beg{enumerate} \item[$(1)$] $P_0^*\gg=\gg$;
  \item[$(2)$] The time-dependent operator  $L_{t, P_t^*\gg}$ 
 generates a unique inhomogenous jump process $\{P_{s,t}^\gg(i,\cdot)\}_{0\le s\le t, i\in \Z_+}$ on $E=\Z_+$ in the sense of Theorem \ref{T0}, such that 
 $$P_t^*\gg= \sum_{i=0}^\infty \gg(i) P_{0,t}^\gg(i,\cdot),\ \ \ t\ge 0. $$ \end{enumerate} \end{defn}
 Let $\P^\gg$ be a distribution dependent  birth-death process generated by $\LL^{a,b}$, and let $X_t$ be the coordinate process on 
 the probability space $(\OO,\scr F,\P):= (D, \B(D),\P^\gg)$, where $\B(D)$ is the Borel $\si$-filed on $D$. Then it is standard to conclude that for any $t\ge 0$ and $i,j\in \Z_+$ with $\P(X_t=i)>0$,
 $$
\lim_{\vv\downarrow 0}\ff 1\vv  \mathbb{P}(X_{t+\vv}=j|X_t=i)=\begin{cases}
	b_t(i, \mathscr{L}_{X_t}),   &  \text{if } j=i+1,\\
	a_t(i, \mathscr{L}_{X_t}),   &  \text{if } j=i-1,\\
	- a_t(i, \mathscr{L}_{X_t})-b_t(i, \mathscr{L}_{X_t}), & \text{if } j=i,\\
	0, & \text{otherwise,}\end{cases} 
$$ where $\L_{X_t}= P_t^*\gg$ is the distribution of $X_t$.

  In applications, we study the  $\LL^{a,b}$-process for distributions in a specific sub-space $\hat{\scr P}$ of $\scr P,$  for instance $\hat{\scr P}=\scr P_p$ for some $p\ge 1$, which leads to the following definition.
 
  \beg{defn}\label{def} We call the $\LL^{a,b}$-problem well-posed for distributions in $\hat{\scr P}$,  if for any $\gg\in \hat{\scr P}$, there is a unique distribution dependent  birth-death process generated by $\LL^{a,b}$
  such that  $P_t^*\gg\in \hat{\scr P}$ for $t\ge 0.$ \end{defn}

The remainder of the paper is organized as follows. In Section \ref{2}, we study the  well-posedness and  exponential ergodicity  of  the $\LL^{a,b}$-problem for distributions in $ \scr P_1$. 
In Section 3, we estimate  the $\W_p$-Lipschitz constant of $P_t^*$ for $p\in [1,\infty)$, which are then applied in Section 4 to characterize   the propagation of chaos for the associated   mean field birth-death systems.  
Finally,  a  result on the well-posedness  of inhomogenous jump processes is presented in  Section 5 as   appendix, which is used in Section 2 to construct the $\LL^{a,b}$-process and its coupling processes.

\section{Well-posedness  and exponential ergodicity  in $\W_1$}\label{2}

 Let $\dd_i$ be the Dirac measure at $i\in \Z_+$, and let  $C^w([0,\infty);\scr P_1)$ be the set of all weakly continuous maps $\mu: [0,\infty)\mapsto \scr P_1$ with
 $$\sup_{t\in [0,T]} \|\mu_t\|_1<\infty,\ \ \ T\in (0,\infty).$$ Let 
$${\rm sgn}(r):=\beg{cases} 1, &\text{if}\ r>0,\\
 -1, &\text{if}\ r<0,\\
 0, &\text{if}\ r=0.\end{cases}$$  
To prove the existence and uniqueness of the $\LL^{a,b}$-process,  we make the following assumption. 
  
\beg{enumerate} \item[$(H_1)$] For any  $i\in \Z_+ $ and $\mu\in C^w([0,\infty);\scr P_1)$,  
$ a_t(i,\mu_t)$ and $b_t(i,\mu_t)$ are  continuous in $t\in [0,\infty)$ locally uniform in $i\in\Z_+$. Moreover,   
there exist $ K_1\in C([0,\infty); \R)$ and $K_2 \in C([0,\infty); [0,\infty))$
such that  
\beq\label{A12} \beg{split}&\big[b_t(i,\mu)-b_t(j,\nu)+a_t(j,\nu)-a_t(i,\mu)\big]{\rm sgn}(i-j)\\
&\quad  + \big( |a_t(i,\mu)-a_t(i,\nu)|+  |b_t(i,\mu)-b_t(i,\nu)|\big)1_{\{i=j\}}\\
&\le K_1(t) |i-j|+ K_2(t)\W_1(\mu,\nu) ,\ \ \ t\ge 0,\ i,j\in \Z_+,\ \mu,\nu\in \scr P_1,\end{split}\end{equation} 
\item[$(H_2)$] There exist    $\theta\in (1,\infty)$, $K_3\in C([0,\infty); [0,\infty))$, $c_0>1$, and  an increasing  function $V: \Z_+\to [1,\infty)$  with set $\{x:V(x)\le c\}$ being compact for any $c>0$ such that  
\beq\label{XC}  V(i+1)\le c_0V(i),\ \ \ i\in \Z_+,\end{equation} 
\beq\label{XC2} 
a_t(i,\dd_0)+ b_t(i,\dd_0) \le K_3(t)V(i),\ \   \ t\ge 0,\  i\in \Z_+, \end{equation}
\beq\label{XC3} \beg{split} 
&\big(b_t(i,\dd_0)-a_t(i,\dd_0)\big) \big(V^\theta(i+1)- V^\theta(i)\big) \\
& \quad + a_t(i,\dd_0) \big(V^\theta(i+1)- 2V^\theta(i)+V^\theta(i-1)\big)\\
&   \le K_3(t)V^\theta(i),\ \ \
t\ge 0,\ i\in \Z_+.\end{split} \end{equation} 
\end{enumerate}

The assumption $(H_1)$ is a standard monotone condition, while $(H_2)$ contains   growth and Lyapunov condition on the birth-death rate at $\dd_0$. 
Before moving on, let us present two examples to illustrate the assumption $(H_2)$.

\beg{exa} Assume $(H_1)$.  \beg{enumerate} 
\item[$(a)$] If there exists   $K\in C([0,\infty);(0,\infty))$   such that 
\beq\label{EX1} a_t(i,\dd_0)\le K(t)  i^2,\ \ \ t\ge 0,\ i\in\Z_+,\end{equation} 
then $(H_2)$ with   $\theta=2$ holds for   some  $K_3\in C([0,\infty);(0,\infty))$ and $V(i):=(1+i)^2.$ 
\item[$(b)$] If there exist  constants $q\ge 1, \vv\in (0,1)$  and  $C_1,C_2,K\in C([0,\infty);(0,\infty))$   such that 
\beq\label{EX2}\beg{split}& b_t(i,\dd_0)-a_t(i,\dd_0)\le C_1(t)- C_2(t) i^q,\\
&  a_t(i,\dd_0)\le K(t)  i^{q+\vv},\ \ \ t\ge 0,\ i\in\Z_+,\end{split} \end{equation} 
with $K(t)\le C_2(t)$ for any $t\ge 0$, then $(H_2)$ with $\theta=2$ holds for some  $K_3\in C([0,\infty);(0,\infty))$  and $V(i):=(1+i)^{q+\varepsilon}$.
\end{enumerate} 
\end{exa} 

\beg{proof}  Obviously,   for any $p\in (0,\infty),$ $V(i):=(1+i)^p $   satisfies \eqref{XC} for some constant $c_0>1.$ 
By  \eqref{A12} with $\mu=\nu=\dd_0$ and $j=0$, we obtain    
\beq\label{BB0} b_t(i,\dd_0)- a_t(i,\dd_0)\le b_t(0,\dd_0)+ K_1(t) i \le C_1(t) (1+i),\ \ \ i\in \Z_+,\ t\in [0,T]\end{equation}
for $C_1(t):= b_t(0,\dd_0)+ K_1(t).$   

$(a)$  If \eqref{EX1} holds, then \eqref{BB0} implies the second condition in $(H_2)$   for $K_3(t)= 2K(t)+ C_1(t)$ and $V(i):=(1+i)^2.$ 
Moreover, for any $p\in [1,\infty)$, we find a constant $c_1>1$ such that $V(i):=(1+i)^p$ satisfies 
\beq\label{RE} \beg{split} &  V^2(i+1)- V^2(i) \\
&\quad = (2+i)^{2 p}-(1+i)^{2 p}\le c_1 (1+i)^{2 p-1}\\
&V^2(i+1)- 2V^2(i)+V^2(i-1)\\
&\quad = (2+i)^{2 p}-2(1+i)^{2 p}+i^{2p}\le c_1 (1+i)^{2 p-2},\ \ i\in \Z_+.\end{split} \end{equation} 
Letting $p=2$ and combining this with 
\eqref{EX1} and  \eqref{BB0}, we derive  
 \beg{align*} &\big(b_t(i,\dd_0)-a_t(i,\dd_0)\big) \big(V^2(i+1)- V^2(i)\big) 
 + a_t(i,\dd_0) \big(V^2(i+1)- 2V^2(i)+V^2(i-1)\big)\\
 &\le C_1(t) c_1(1+i)^{2p} + K(t)i^{2} c_1 (1+i)^{2p-2} \le \big(C_1(t)+K(t))c_1V^ 2(i),\ \ t\ge 0,\ i\in \Z_+.\end{align*}
 So, the third condition in $(H_2)$   holds for some $K_3$ and $V(i)=(1+i)^2$. 
 
 $(b)$ Let $V(i)= (1+i)^{q+\vv}$. By the second inequality in \eqref{EX2} and \eqref{BB0}, we obtain   the second condition in $(H_2)$ for some $K_3$. 
Next, by \eqref{BB0} and \eqref{RE} with $p=q+\vv$ for the present $V$, we obtain
 \beg{align*} &\big(b_t(i,\dd_0)-a_t(i,\dd_0)\big) \big(V^2(i+1)- V^2(i)\big) 
 + a_t(i,\dd_0) \big(V^2(i+1)- 2V^2(i)+V^2(i-1)\big)\\
 &\le  C_1(t)c_1 (1+i)^{2 (q+\vv)-1}  - C_2(t)c_1(1+i)^{2 (q+\vv)-1} i^{q}  + K(t)i^{q+\vv} c_1 (1+i)^{2(q+\vv)-2}
 \\&\le C_1(t)c_1 (1+i)^{2 (q+\vv)}
 ,\ \ t\ge 0,\ i\in \Z_+.\end{align*}
 Here we use the fact that $K(t)\le C_2(t)$ and $\varepsilon \in (0, 1)$. This implies the third condition in $(H_2)$ for $\theta=2$ and 
 $V( i)=(1+i)^{q+\vv}.$ 
 
\end{proof} 

\beg{thm}\label{T1} Assume $(H_1)$ and $(H_2)$. Then the following assertions hold.
\beg{enumerate} \item[$(1)$] The $\LL^{a,b}$-problem is well-posed for distributions in $\scr P_1$, and 
\beq\label{ES0} \|P_t^*\mu\|_1 \le \e^{\int_0^t (K_1+K_2)(s)\d s } \|\mu\|_1  + \int_0^t b_s(0,\dd_0)\e^{\int_s^t (K_1+K_2)(r)\d r}   \d s,\ \ t\ge 0,\end{equation} 
 \beq\label{ES1} \W_1(P_t^*\mu, P_t^*\nu)\le\e^{\int_0^t (K_1(s)+K_2(s))\d s} \W_1(\mu,\nu),\ \ \ t\ge 0,\ \mu,\nu\in \scr P_1.\end{equation} 
\item[$(2)$]   If  $(a_t(i,\mu), b_t(i,\mu))=(a(i,\mu),b (i,\mu))$ does not depend on $t$  such that  $(H_1)$ holds for some constants $K_1\in\R$ and $K_2\ge 0$ with $K_1+K_2<0,$ 
then   $P_t^*$ has a unique invariant probability measure $\bar \mu\in \scr P_1$, i.e. $P_t^*\bar \mu=\bar \mu$ for all $t\ge 0$, and 
\beq\label{ES1'} \W_1(P_t^*\mu,  \bar\mu) \le\e^{(K_1+K_2)t} \W_1(\mu,\bar\mu),\ \ \ t\ge 0,\ \mu  \in \scr P_1.\end{equation} 
\end{enumerate} 
  \end{thm}
 
 \beg{proof}  (1) It suffices to consider the $\LL^{a,b}$-problem  up to an arbitrarily  finite  time. So, in the following, we fix $T\in (0,\infty)$ and an initial distribution $\mu\in \scr P_1$.  
 
 (a) For  $\mu\in \scr P_1$ and $T\in (0,\infty)$, let
 $$\C_T^\mu:= \bigg\{\gg\in C^w([0,T];\scr P):\ \gg_0=\mu,\ H_T(\gg):=\sup_{t\in [0,T]} \|\gg_t\|_1<\infty\bigg\}.$$
 For fixed $\gg\in \C_T^\mu$, \eqref{A12} with $\mu=\nu=\gg_t$ and $j=0$  implies 
\begin{align}\label{BB} b_t(i,\gg_t)- a_t(i,\gg_t)&\le b_t(0,\gg_t)+ K_1(t) i\notag\\&\leq b_t\left(0, \delta_0\right)+K_2(t)\left\|\gamma_t\right\|_1+K_1(t) i \leq c_1 +c_1 i,\ \ \ i\in \Z_+,\ t\in [0,T]\end{align}
for some constant $c_1>0$, since $b_t(0,\gg_t)$ and $ K_1(t)$ are bounded in $t\in [0,T].$  
Let
\beq\label{LT} L_t^\gg f(i):= a_t(i,\gg_t)(f(i-1)- f(i)) +b_t(i,\gg_t)(f(i+1)-f(i)),\ \ f\in \B_b(\Z_+).\end{equation} 
 Next, taking $i=j$ in \eqref{A12} we derive
\beq\label{*1} |a_t(i,\mu_1)-a_t(i,\mu_2)| +  |b_t(i,\mu_1)-b_t(i,\mu_2)| \le K_2(t) \W_1(\mu_1,\mu_2),\ \   \mu_1,\mu_2\in \scr P_1.\end{equation}
This together with $\eqref{XC2}$   yields 
\beq\label{QO} q_t(i):= a_t(i,\gg_t)+ b_t(i,\gg_t)\le c_2(\gg)V(i),\ \ \  \ i\in\Z_+,\ t\in [0,T]\end{equation} 
for some constant $c_2(\gg)$ depending on $\gg$. \newpage

Next, by \eqref{XC}, $\eqref{XC3}$ and \eqref{*1}, we find  a constant  $c_3(\gg)>0$ such that 
\beg{align}\label{UU} &L_t^\gg V^\theta (i) = L_t^{\dd_0}V^\theta (i)+ (a_t(i,\gg_t)-a_t(i,\dd_0)) (V^\theta (i-1)-V^\theta (i)) \notag\\
&\qquad \qquad\qquad + (b_t(i,\gg_t)-b_t(i,\dd_0)) (V^\theta (i+1)-V^\theta (i))\\
&\le K_3(t) V^\theta(i)+ (c_0^{\theta}-1)K_2(t)\|\gg_t\|_1   V^\theta (i)\le  c_3(\gg) V^\theta (i),\ \ t\in [0,T],\ i\in\Z_+.\notag\end{align}
 So,  by Theorem \ref{T0}, there is a unique inhomogenous birth-death process  $$\{P_{s,t}^\gg(i,\cdot)\}_{i\in \Z_+, 0\le s\le t\le T}$$  generated by
$\{L_t^\gg\}_{t\in [0,T]}.$ 
  Let
$$(\Phi \gg)_t=\Phi_t \gg := \sum_{i=0}^\infty \mu(i) P_{s,t}^\gg(i,\cdot),\ \ \ t\in [0,T].$$
If we can verify that $\Phi: \C_T^\mu\to \C_T^\mu$ has a unique fixed point $\bar\gg$, then   the inhomogenous birth-death process generated by $\{L_t^{\bar\gg}\}_{t\in [0,T]}$ is the 
unique $\LL^{a,b}$-process up to time $T$ with initial distribution $\mu$. 

(b) We first show $\Phi\gg\in \C_T^\mu$ for $\gg\in \C_T^\mu$.  The weak continuity of $t\mapsto (\Phi\gg)_t$ is  standard  for the birth-death process. 
By \eqref{BB} we have 
$$\sum_{j=0}^\infty (1+j ) P_{s,t}^\gg(i,\{j\}) \le (1+i) \e^{c_1(t-s) },\ \ \ 0\le s\le t\le T,\ i\in\Z_+,$$
so that  
\beg{align*} \|\Phi_t\gg\|_1=\sum_{i,j=0}^\infty \mu(i) j P_{s,t}^\gg(i,\{j\}) \le   (1+ \|\mu\|_1) \e^{c_1T}<\infty,\ \ \ t\in [0,T].\end{align*}
 Then $ \Phi \gg\in \C_T^\mu$. 

To see that $\Phi$ has a unique fixed point, we need only to show that $\Phi$ is contractive under the 
complete metric 
$$\rr_\ll(\gg,\tt\gg):=\sup_{t\in [0,T]} \e^{-\ll t} \W_1(\gg_t,\tt\gg_t),\ \ \gg,\tt\gg\in \C_T^\mu$$
for some constant $\ll>0$.  To this end, we consider the following synchronized coupling operator for $L_t^\gg$ and $L_t^{\tt\gg}$: for  $t\ge 0, i,j\in \Z_+$ and $f\in \B(\Z_+^2),$  
\beq\label{CPO}\beg{split} L_t^{\gg,\tt\gg} f(i,j)  := &\, \big[b_t(i,\gg_t)\land b_t(j,\tt\gg_t)\big] \big(f(i+1,j+1)-f(i,j)\big) \\
&+ \big[a_t(i,\gg_t)\land a_t(j,\tt\gg_t)\big] \big(f(i-1,j-1)-f(i,j)\big)\\ 
&+ \big[b_t(i,\gg_t)- b_t(j,\tt\gg_t)\big]^+ \big(f(i+1,j)-f(i,j)\big) \\
&+ \big[b_t(i,\gg_t)- b_t(j,\tt\gg_t)\big]^- \big(f(i,j+1)-f(i,j)\big) \\
& +\big[a_t(i,\gg_t)- a_t(j,\tt\gg_t)\big]^+ \big(f(i-1,j)-f(i,j)\big) \\
&+ \big[a_t(i,\gg_t)- a_t(j,\tt\gg_t)\big]^- \big(f(i,j-1)-f(i,j)\big).\end{split}\end{equation} 
By taking $$\bar V(i,j) := \big(V(i)^\theta+V(j)^\theta\big)^{\ff 1 \theta},$$
the marginal property of the coupling operator and \eqref{UU} imply
$$L_t^{\gg,\tt\gg} \bar V^\theta (i,j)\le (c_3(\gg)+c_3(\tt\gg))\bar V^\theta (i,j),\ \ \ t\in [0,T],\ i,j\in \Z_+.$$
Moreover, by \eqref{QO},    $q$-pair of  the coupling operator satisfies
$$\bar q_t(i,j):= a_t(i,\gg_t)+ b_t(i,\gg_t)+ a_t(j,\tt \gg_t)+ b_t(j,\tt \gg_t)\le \big(c_2(\gg)+c_2(\tt\gg) \big)\bar V(i,j).$$
So, by Theorem \ref{T0},   there exists a unique  coupling process up to time $T$  generated by  $\{L_t^{\gg,\tt\gg}\}_{t\in [0,T]}$. 
  Let 
  $ (X_t^\gg, Y_t^{\tt\gg} )$ be the coordinate process on $(\OO,\F,\P)$, where 
  $$\OO:= D([0,T];\Z_+\times\Z_+)$$ is the space of  c$\grave{a}$dl$\grave{a}$g paths from $[0,T]$ to $\Z_+\times\Z_+,$  $\F$ is the Borel $\si$-field of $\OO$, and $\P$
  is the distribution of the coupling process with initial distribution 
\beq\label{S1} \pi_{ij} := \mu(i )1_{\{i=j\}},\ \ i,j\in \Z_+.\end{equation}  By \eqref{A12}, for $\rr (i,j):=|i-j|$  we obtain 
\beq\label{CP} \beg{split} L_t^{\gg,\tt\gg}\rr (i,j)&=  
  \big[b_t(i,\gg_t)- b_t(j,\tt\gg_t)\big]^+ \big(\rr(i+1,j)-\rr(i,j)\big)\\
 &\quad  + \big[b_t(i,\gg_t)- b_t(j,\tt\gg_t)\big]^- \big(\rr(i,j+1)-\rr(i,j)\big) \\
&\quad +  \big[a_t(i,\gg_t)- a_t(j,\tt\gg_t)\big]^+ \big(\rr(i-1,j)-\rr(i,j)\big)\\
 &\quad  + \big[a_t(i,\gg_t)- a_t(j,\tt\gg_t)\big]^- \big(\rr(i,j-1)-\rr(i,j)\big)\\
&={\rm sgn}(i-j)  \big[ b_t(i,\gg_t)- b_t(j,\tt\gg_t) -a_t(i,\gg_t)+ a_t(j,\tt\gg_t)\big] \\
&\quad + 1_{\{i=j\}}   \big[|b_t(i,\gg_t)- b_t(i,\tt\gg_t)|+ |a_t(i,\gg_t)- a_t(i,\tt\gg_t)| \big]\\
&\le K_1(t)|i-j|+ K_2(t) \W_1(\gg_t,\tt\gg_t),\ \ \ t\ge 0,\ i,j\in \Z_+. \end{split}\end{equation} 
This implies 
$$\E[|X_t^\gg-Y_t^{\tt\gg}|] -\E[|X_s^\gg-Y_s^{\tt\gg}|] \le \int_s^t \Big(K_1(r) \E[|X_r^\gg-Y_r^{\tt\gg}|] +K_2(r) \W_1(\gg_r,\tt\gg_r)\Big)\d r,\ \ t\ge s\ge 0.$$
By  Gronwall's inequality and \eqref{S1},   we derive
$$ \W_1(\Phi_t\gg, \Phi_t\tt\gg ) \le \E[|X_t^\gg-Y_t^{\tt\gg}|] \le \int_0^t \e^{\int_0^t K_1(r)\d r} K_2(r) \W_1(\gg_r,\tt\gg_r)\d r,\ \ t\in [0,T].$$
So, when $\ll>0$ is large enough,
\beg{align*} &\rr_\ll(\Phi\gg,\Phi\tt\gg):=\sup_{t\in [0,T]} \e^{-\ll t}  \W_1(\Phi_t\gg, \Phi_t\tt\gg )  \\
&\le \Big(\sup_{t\in [0,T]} K_2(t)\Big) \e^{\int_0^T |K_1(r)|\d r} \rr_\ll(\gg,\tt\gg)    \int_0^T  \e^{-\ll r}\d r\\
&\le \ff 1 2  \rr_\ll(\gg,\tt\gg),\ \ \ \gg,\tt\gg\in \C_T^\mu.\end{align*} 
 Thus,  when $\ll>0$ is large enough, $\Phi$ is contractive in $\rr_\ll,$ so that it has a unique fixed point in $\C_T^\mu$. Since $T\in (0,\infty)$ is arbitrary,  there is a unique $\LL^{a,b}$-process with initial distribution $\mu$.

(c)  To prove \eqref{ES0}, let $\gg_t=P_t^*\mu$, so that \eqref{BB} implies 
 $$\|\gg_t\|_1 -\|\gg_s\|_1  \le \int_s^t \big[b_r(0,\dd_0)+ (K_1+K_2)(r) \| \gg_r\|_1\big]\d r,\ \ 0\le s\le t<\infty.$$
By Gronwall's inequality, this implies \eqref{ES0}.

Let $\gg_t=P_t^*\mu, \tt \gg_t= P_t^*\nu$, and let $(X_t^\gg,Y_t^{\tt\gg})$ be the coupling process constructed as in above, but with initial distribution  $\pi=(\pi_{ij})_{i,j\in \Z_+} \in \C(\mu,\nu)$   such that
 \beq\label{SS} \W_1(\mu,\nu)=\sum_{i,j=0}^\infty \pi_{ij} |i-j|=\E[|X_0^\gg-Y_0^{\tt\gg}|].\end{equation} 
On the other hand, we have 
$$  \W_1(P_t^*\mu,P_t^*\nu)\le \E [|X_t^\gg-Y_t^{\tt\gg}|] =:h_t.$$ Then  
\eqref{CP}  implies 
$$h_t -h_s\le \int_s^t \big[K_1(r) h_r + K_2(r) \W_1(\gg_r,\tt\gg_r)\big]\d r\le \int_s^t (K_1+K_2)(r)h_r\d r,\ \ 0\le s\le t<\infty.$$
Combining this with  \eqref{SS} and Gronwall's inequality, we derive \eqref{ES1}.

(2)  When $(a_t(i,\mu), b_t(i,\mu))=(a(i,\mu),b(i,\mu))$ such that $K_1$ and $K_2$ are constant with $\ll:= -(K_1+K_2)>0$, then
\eqref{ES1} becomes the exponential contraction of $P_t^*$: 
$$  \W_1(P_t^*\mu,P_t^*\nu)\le \e^{-\ll t}\W_1(\mu,\nu),\ \ t\ge 0, \mu,\nu\in \scr P_1.$$ 
 According to \cite[Theorem 5.1.1]{WR},   this together with \eqref{ES0} implies that $P_t^*$ has a unique invariant probability measure $\bar \mu\in \scr P_1$ and \eqref{ES1'} holds.  
 \end{proof} 
 
 \paragraph{Remark 2.1.} In the context of  DDSDEs,  the exponential ergodicity  has been proved for the partially dissipative case and non-dissipative case, see \cite{wfy23ber, wfy23spa}. 
 In the same spirit, it should be  possibly to derive the exponential ergodicity for the present model by replacing \eqref{A12}   with $K_1+K_2<0$ by the weaker condition 
\beq\label{A120} \beg{split}&\big[b(i,\mu)-b(j,\nu)+a(j,\nu)-a(i,\mu)\big]{\rm sgn}(i-j)\\
&\quad  + \big( |a(i,\mu)-a(i,\nu)|+  |b(i,\mu)-b(i,\nu)|\big)1_{\{i=j\}}\\
&\le K+ K_1 |i-j|1_{\{|i-j|\ge N\}} + K_2  \W_1(\mu,\nu) ,\ \ \ t\ge 0,\ i,j\in \Z_+,\ \mu,\nu\in \scr P_1,\end{split}\end{equation} 
for some constants  $K>0, N\in \N, K_1<0 $ and small enough $K_2>0$.     
To this end, one should adopt the coupling by reflection rather than the synchronized coupling in \eqref{CPO}.

 \section{Lipschitz  continuity in $\W_p$ }
 In this section, for fixed $p\in (1,\infty)$, we study the   the Lipschitz continuity of  the maps 
\beg{align*} &P_t^*: \scr P_p\to\scr P_p,\\ 
&\scr P_p\ni\mu\mapsto P_tf(\mu) := \int_{\Z_+} f\d (P_t^*\mu),\ \ \ f\in \B_b(\Z_+) \end{align*}  for $t> 0.$ 
 To this end, we make  the following assumption.
 \beg{enumerate}   
 \item[$(H_3)$] Let $p\in (1,\infty)$. There exist  continuous functions  $\bb,\beta_1\in C([0,\infty); (0,\infty)), \bb_k \in C([0,\infty); \R), 2\le k\le 3,$ such that  
   \beg{align*}& p\big[(i+1)^{p-1}b_t(i,\mu)-(i-1)^{p-1} a_t(i,\mu)\big] 
  \le \bb_1(t) + \bb_2(t) i^p+\bb_3(t) \|\mu\|_1^p,\\
 &|a_t(i,\mu)-a_t(j,\nu)|+ |b_t(i,\mu)-b_t(j,\nu)|\le \bb(t)\big(|i-j|+ \W_1(\mu,\nu)\big),\\
 &\  \ i,j\in\Z_+,\ \ t \ge 0,\ \  \mu,\nu\in \scr P_1.\end{align*}
  \end{enumerate}
 
 We note that if   $(H_1)$ and one of  \eqref{EX1} and \eqref{EX2} hold, then     $(H_2)$ and $(H_3)$ are satisfied, where $(H_2)$  is already verified in Example 2.1 and $(H_3)$ can be deduced similarly.

 \beg{thm}\label{T2} Let $p\in (1,\infty)$ and assume $(H_1)$-$(H_3)$. 
 Then  
    \beq\label{ES20}  \|P_t^*\mu\|_p^p \le \e^{ \int_0^t(\bb_2+\bb_3)(s)\d s}\|\mu\|_p^p + \int_0^t \bb_1(s)\e^{ \int_s^t(\bb_2+\bb_3)(r)\d r} \d s,\ \ t\ge 0,\ \mu\in \scr P_p.  \end{equation}
  \beq\label{ES2} \W_p(P_t^*\mu, P_t^*\nu)\le \e^{ 2^p\int_0^t  \bb (s) \d s} \W_p(\mu,\nu),\ \ \ t\ge 0,\ \mu,\nu\in \scr P_p.\end{equation}  
Consequently,  letting
$$|\nn f(i)|:=\sup_{j\ne i} \ff{|f(j)-f(i)|}{|i-j|},\ \ i\in \Z_+,\ f\in \B_b(\Z_+),$$
for any $p>1$, we have 
  \beq\label{ES2'}  \sup_{\nu\ne \mu} \ff{|P_tf(\mu)- P_tf(\nu)|}{\W_p (\mu,\nu)} \le \e^{2^{p}\int_0^t  \bb (s) \d s}\big(P_t|\nn f|^{\ff p{p-1}}(\mu)\big)^{\ff{p-1}p},\ \ t\ge 0,\ \mu\in \scr P_p.\end{equation} 
\end{thm}

 \beg{proof} Let $\mu,\nu\in \scr P_1$ and $\gg_t=P_t^*\mu, \tt\gg_t=P_t^*\nu,\ t\ge 0.$  For fixed constant $p>1$,
 let $g_p(i):= i^p, i\in\Z_+$. By $(H3)$, the operator $L_t^\gg$ defined in \eqref{LT} satisfies
 \beg{align*} &L_t^\gg g_p(i)= b_t(i,\gg_t)\big[(i+1)^p-i^p\big]+ a_t(i,\gg_t)\big[(i-1)^p-i^p\big]  \\
 & \le p\big[(i+1)^{p-1}b_t(i,\gg_t)-(i-1)^{p-1} a_t(i,\gg_t)\big]\\
  &\le \bb_1(t) +\bb_2(t) i^p+\bb_3(t) \|\gg_t\|_1^p,\ \ t\ge 0,\ i\in\Z_+.\end{align*}
This and \eqref{ES0} imply that $\sup_{t\in [0,T]} \|\gg_t\|_p<\infty$ for $T\in (0,\infty)$ and 
 $$\|\gg_t\|_p^p \le \|\mu\|_p^p +\int_0^t   \big[\bb_1(s)+(\bb_2+\bb_3)(s) \|\gg_s\|_p^p\big]\d s,\ \ t\ge 0,$$
 which implies \eqref{ES20} by Gronwall's inequality. 
 
 Next, let $f(i,j)=|i-j|^p$. 
   By   $(H_3)$ and 
\beg{align*}&\max\Big\{ \big||i-j+1|^p-|i-j|^p\big|,\ \big||i-j-1|^p-|i-j|^p\big|\Big\} \\
&\le p (|i-j|+1)^{p-1} \le p 2^{p-1} (1\lor |i-j|)^{p-1},\end{align*} 
 the coupling operator $L_t^{\gg,\tt\gg}$ defined in \eqref{CPO} satisfies 
\beq\label{II}   \beg{split}  L_t^{\gg,\tt\gg}f (i,j)&=  
  \big[b_t(i,\gg_t)- b_t(j,\tt\gg_t)\big]^+ \big(|i-j+1|^p-|i-j|^p  \big)\\
 &\quad  + \big[b_t(i,\gg_t)- b_t(j,\tt\gg_t)\big]^- \big( |i-j-1|^p- |i-j|^p\big) \\
&\quad +  \big[a_t(i,\gg_t)- a_t(j,\tt\gg_t)\big]^+ \big(|i-j-1|^p-|i-j|^p\big)\\
 &\quad  + \big[a_t(i,\gg_t)- a_t(j,\tt\gg_t)\big]^- \big(|i-j+1|^p-|i-j|^p\big)\\
&\le p 2^{p-1} (1\lor |i-j|)^{p-1}  \big(|b_t(i,\gg_t)- b_t(j,\tt\gg_t)|+ |a_t(i,\gg_t)-a_t(j,\tt\gg_t)|\big) \\ 
&\le p 2^{p-1}   \bb(t)   (1\lor |i-j|)^{p-1} \big(|i-j|+\W_1(\gg_t,\tt\gg_t)\big),\ \ \ t\ge 0,\ i,j\in \Z_+. \end{split} \end{equation} 
By Theorem \ref{T0}, let  $(X_t^\gg, Y_t^{\tt\gg})$ be generated by $L_t^{\gg,\tt\gg}$ with initial distribution $\pi\in \C(\mu,\nu)$ such that
 $$\W_p(\mu,\nu)^p= \sum_{i,j=0}^\infty \pi_{ij}  |i-j|^p=\E[|X_0^\gg-Y_0^{\tt\gg}|^p].$$
Noting that  
$$|i-j|(1\lor |i-j|)^{p-1}= |i-j|^p,\ \ \ i,j\in \Z_+,$$
which together with   $\W_1(\gg_t,\tt\gg_t)\le \E[|X_t^\gg-Y_t^{\tt\gg}|]$ and the FKG inequality implies
\beg{align*}&\E[(1\lor |X_t^\gg-Y_t^{\tt\gg}|)^{p-1} ] \W_1(\gg_t,\tt\gg_t)\le \E[(1\lor |X_t^\gg-Y_t^{\tt\gg}|)^{p-1} ] \E[|X_t^\gg-Y_t^{\tt\gg}|]\\
& \le \E[(1\lor |X_t^\gg-Y_t^{\tt\gg}|)^{p-1}  |X_t^\gg-Y_t^{\tt\gg}|] 
= \E[|X_t^\gg-Y_t^{\tt\gg}|^p].\end{align*} 
We deduce from \eqref{II} that 
\beg{align*}  &\E[|X_t^\gg-Y_t^{\tt\gg}|^p]\le  \W_p(\mu,\nu)^p\\
&+  p 2^{p-1}  
  \E \int_0^t   \bb(s)  \big(1\lor |X_s^\gg-Y_s^{\tt\gg}|\big)^{p-1} \big(|X_s^\gg-Y_s^{\tt\gg}|+\W_1(\gg_s,\tt\gg_s)\big) \d s\\
  &\le   \W_p(\mu,\nu)^p+  p 2^p   
  \int_0^t   \bb(s) \E[ |X_s^\gg-Y_s^{\tt\gg}|^{p} ] \d s,\ \ t\ge 0.
 \end{align*} 
   Since  $\E[|X_t^\gg-Y_t^{\tt\gg}|^p]<\infty$ is ensured by \eqref{ES20},  by Gronwall's inequality we obtain 
 \beq\label{ES2''}\Big( \E[|X_t^\gg-Y_t^{\tt\gg}|^p]\Big)^{\ff 1 p} \le \e^{  2^p \int_0^t  \bb (s) \d s} \W_p(\mu,\nu).\end{equation}
 This implies \eqref{ES2}  since 
  $  \W_p(P_t^*\mu, P_t^*\nu)^p \le \E[|X_t^\gg-Y_t^{\tt\gg}|^p].$  
      
    Finally, by   H\"older's inequality, we obtain  
    \beg{align*}  &\ff{|P_tf(\mu)- P_tf(\nu)|}{\W_p (\mu,\nu)}  \le \ff{\E[|f(X_t^\gg)-f(Y_t^{\tt\gg})|]}  {\W_p (\mu,\nu)} \\
   & \le \E\bigg[|\nn f|(X_t^\gg)\cdot\ff{|X_t^\gg-Y_t^{\tt \gg }|}{  \W_p (\mu,\nu)} \bigg]\\
   &\le \big(P_t |\nn f|^{\ff p{p-1}}\big)^{\ff{p-1}p} \ff{(\E[|X_t^\gg-Y_t^{\tt \gg}|^p])^{\ff 1 p}}{ \W_p (\mu,\nu)}.\end{align*}
   Combining this with     \eqref{ES2''}  we derive \eqref{ES2'} .
  \end{proof} 

 
\section{Propagation of chaos}\label{poc}
 
Let $N\ge 2$,  and $\scr P_1(N)$ be the space of all probability measures $\nu$ on $\Z_+^N$ with  $\nu(\rr_N(x,o))<\infty$, where $o=(0,0,\cdots,0)$ and 
$$\rr_N(x,y)=|x-y|:=\sum_{l=1}^N|x_l-y_l|,\ \ x,y\in \Z_+^N.$$

A measure $\nu\in \scr P_1(N)$ is called symmetric, if 
 for any permutation $l:=(l_1,\cdots,l_N)$ of $(1,\cdots, N)$,
$\nu$ is invariant under the map
$$\Z_+^N\ni x=(x_1,\cdots,x_N)\mapsto \pi_l x:= (x_{l_1},\cdots, x_{l_N})\in \Z_+^N.$$
For any $x\in \Z_+^N$, let
$$\mu^N (x):=  \ff 1 N \sum_{l=1}^N \dd_{x_l},$$
where $\dd_{x_l}$ is the Dirac measure at $x_l\in\Z_+$. Let $e_l\in \Z_+^N$ with the $l$-th component $1$ and others $0$.

We consider the mean field particle systems $X_t^N=(X_t^{N,1},\cdots, X_t^{N,N})$ on $\Z_+^N$  generated by
$$ L_t^Nf(x)=\sum_{l=1}^N \Big[b_t\big({x_l}, \mu^N(x)\big)\big(f(x+e_l)-f(x)\big)+a_t\big({x_l},\mu^N(x)\big)\big(f(x-e_l)-f(x)\big)\Big]$$ for $
  f\in \B(\Z_+^N),\ x\in \Z_+^N. $ 
For any $1\le k\le N,$ let  
  $$P_t^{N*k}\nu =\L_{(X_t^{N,1},\cdots, X_t^{N,k})}\ \text{for}\ \L_{X_0^N}=\nu,\ \ t\ge 0,\ \nu\in \scr P_1(N).$$
  When $k=N$ we simply denote $P_t^{N*}= P_t^{N*N}$.
  
For any $\mu_0\in \scr P_2$,   let $\mu_0^{\otimes N}$ be the $N$ times product measure of $\mu_0$.
Under $(H_1)$-$(H3)$ for $p=2$, let 
\beg{align*} & h_t(\mu_0):=\bigg( \e^{ \int_0^t(\bb_2+\bb_3)(s)\d s}\|\mu_0\|_2^2  + \int_0^t \bb_1(s)\e^{ \int_s^t(\bb_2+\bb_3)(r)\d r} \d s\bigg)^{\ff 1 2},\\
&H_t(\mu_0):=   1+h_t(\mu_0)+  \int_0^t \big(1+h_s(\mu_0)\big) K_2(s) \e^{\int_s^t K_1(r)\d r}  \d s,\ \ t\ge 0.\end{align*} 
  Let $\W_1$ be the $1$-Wasserstein distance on $\scr P_1(N)$ induced by $\rr_N$. 
 
\beg{thm}\label{T3} Assume $(H_1)$-$(H_3)$   for $p=2$. 
 \begin{enumerate}
 	\item[$(1)$] There exists a constant $c>0$ such that for any $N\ge 2$, $\mu_0\in \scr P_2$ and $\nu\in \scr P_1(N)$, we have 
 	\beq\label{FN1}\beg{split}& \W_1\big(P_t^{N*}\nu, (P_t^*\mu_0)^{\otimes N}\big)\\
 	&\le \e^{\int_0^t (K_1+K_2)(s)\d s} \W_1(\nu,\mu_0^{\otimes N}) 
 	+ c\ss N \int_0^t \e^{\int_s^t K_1(r)\d r} K_2(s) H_s(\mu_0)\d s.\end{split} \end{equation}
 Consequently, if $\nu$ is symmetric, then for any $N\ge 2$ and $ 1\le k\le N$, 
\beq\label{FN2}\beg{split}& \W_1\big(P_t^{N*k}\nu, (P_t^*\mu_0)^{\otimes k}\big)\\
&\le \ff k N \e^{\int_0^t (K_1+K_2)(s)\d s} \W_1(\nu,\mu_0^{\otimes N}) 
+ \ff{ck}{\ss N} \int_0^t \e^{\int_s^t K_1(r)\d r} K_2(s) H_s(\mu_0)\d s.\end{split}\end{equation} 
    \item[$(2)$] If  $(a_t(i,\mu), b_t(i,\mu))=(a(i,\mu),b (i,\mu))$ does not depend on $t$  such that  $(H_1)$ holds with $K_1+K_2<0$ and $\beta_2+\beta_3<0$, then for any $N\ge 2$, $1\le k\le N$, $\mu_0\in\mathscr{P}_2$ and $\nu\in\mathscr{P}_1(N)$ being symmetric, there exists a constant $c>0$ such that
    	\beq\label{FN3} \W_1\big(P_t^{N*k}\nu, (P_t^*\mu_0)^{\otimes k}\big)
    \le \ff k N \e^{ (K_1+K_2)t} \W_1(\nu,\mu_0^{\otimes N}) 
    + c\ff{k}{\ss N} (1+\|\mu_0\|_2). \end{equation}
 \end{enumerate}
 
\end{thm} 
   
   \beg{proof} It is standard that \eqref{FN1} implies \eqref{FN2}. By the triangle inequality,
   we have
   $$\W_1\big(P_t^{N*}\nu, (P_t^*\mu_0)^{\otimes N}\big)\le \W_1\big(P_t^{N*}\nu, P_t^{N*}\mu_0^{\otimes N} \big) + \W_1\big(P_t^{N*}\mu_0^{\otimes N},(P_t^*\mu_0)^{\otimes N}\big).$$
Then \eqref{FN1} follows from \eqref{AB2} and \eqref{Y1} for $\tt\nu=\mu_0^{\otimes N}$, which are proved in the following two lemmas. \eqref{FN3} is easily followed by \eqref{FN2} by direct calculations.
   
   \end{proof} 
 
\beg{lem} \label{LN2}  Under the assumption of Theorem \ref{T3}, for any $\mu_0\in \scr P_2$, let $\mu_t:=P_t^*\mu_0$ for $t\ge 0$, and 
 moreover, for any $N\ge 2$, let   $X_t^N$ have    initial distribution $\mu_0^{\otimes N}$.  Then there exists a constant $c>0$ such that for any $t\ge 0$ and $X_t^N$ having initial distribution $\mu_0^{\otimes N},$
we have 
\beq\label{AB1} \E\big[\W_1(\mu^N(X_t^N), \mu_t)\big] \le \ff{c}{\ss N}  H_t(\mu_0), \end{equation}
\beq\label{AB2}  \W_1\big(P_t^{N*}\mu_0^{\otimes N}, (P_t^*\mu_0)^{\otimes N}\big)\le  c\ss N  \int_0^t \e^{\int_s^t K_1(r)\d r} K_2(s) H_s(\mu_0) \d s.\end{equation} 
 
\end{lem}

\beg{proof}    
Given  the initial value $X_0^N=Y_0^N$ with $\L_{X_0^N}=\mu_0^{\otimes N}$, by Theorem \ref{T0} there is a unique inhomogenous Markov process   $(X_t^N, Y_t^N)$  on $\Z_+^N\times \Z_+^N$  generated by 
$\{\tt L_t^{\mu,N}: t\ge 0\}$, where  for $t\ge 0, f\in \B(\Z_+^N)$ and $x=(x_1,\cdots,x_N), y=(y_1,\cdots,y_N)\in \Z_+^N$,
\beg{align*}   \tt L_t^{\mu,N} f(x,y) :=  \sum_{l=1}^N   &\Big[ \big(b_t(x_l,\mu^N(x))\land b_t(y_l,\mu_t)\big) \big(f(x+e_l,y+e_l)-f(x,y)\big) \\
&+ \big(a_t(x_l,\mu^N(x))  \land a_t(y_l,\mu_t) \big) \big(f(x-e_l,y-e_l)-f(x,y)\big)\\ 
&+  \big(b_t(x_l,\mu^N(x))- b_t(y_l,\mu_t)\big)^+ \big(f(x+e_l,y)-f(x,y)\big)  \\
&+ \big(b_t(x_l,\mu^N(x))- b_t(y_l,\mu_t)\big)^- \big(f(x,y+e_l)-f(x,y)\big) \\
& +\big(a_t(x_l,\mu^N(x))- a_t(y_l,\mu_t)\big)^+ \big(f(x-e_l,y)-f(x,y)\big) \\
&+ \big(a_t(x_l,\mu^N(x))- a_t(y_l,\mu_t)\big)^- \big(f(x,y-e_l)-f(x,y)\big)\Big].\end{align*}
Then 
$$\L_{X_t^N}= P_t^{N*}\mu_0^{\otimes N},\ \ \ \L_{Y_t^N}=  (P_t^*\mu_0)^{\otimes N},\ \ t\ge 0.$$
 Then as in \eqref{CP}, we have 
 \beg{align*}&\tt L_t^{\mu,N} \rr_N (x,y) \le \sum_{l=1}^N\big[K_1(t)|x_l-y_l|+ K_2(t)   \W_1(\mu^N(x),\mu_t)\big]\\
&= K_1(t)\rr_N(x,y)+ N K_2(t)   \W_1(\mu^N(x),\mu_t),\ \ t\ge 0,\ x,y\in \Z_+^N.\end{align*}
Noting that $X_0^N=Y_0^N$, this implies
\beq\label{K0}   \E[|X_t^N-Y_t^N|] \le  N\int_0^t \e^{\int_s^t K_1(r)\d r} K_2(s)   \E[ \W_1(\mu^N(X_s^N),\mu_s)]\d s.\end{equation} 
Consequently, 
 \beq\label{K1}\beg{split} &\E[ \W_1(\mu^N(X_t^N), \mu^N(Y_t^N))] \le \ff 1 N  \E[|X_t^N-Y_t^N|]  \\
&\le  \int_0^t \e^{\int_s^t K_1(r)\d r} K_2(s) \E[\W_1(\mu^N(X_s^N),\mu_s)] \d s,\ \ t\ge 0.\end{split}\end{equation} 
 On the other hand, noting that the components $ \{Y_t^{i}\}_{1\le i\le N} $ of $Y_t^N$ are i.i.d. with distribution $\mu_t$, by \cite[Lemma 2.1]{huang2024}, there exists a constant $c>0$ such that 
 $$\E[\W_1(\mu^N(Y_t^N),\mu_t)]\le  \ff c{\ss N} \big(1+\|\mu_t\|_2\big),\ \ t\ge 0, \ N\ge 2.$$
 Combining this with \eqref{FN1}  and the triangle inequality, we obtain
\beg{align*} &\E[\W_1(\mu^N(X_t^N), \mu_t)]\le  \E[ \W_1(\mu^N(X_t^N), \mu^N(Y_t^N))] + \E[\W_1(\mu^N(Y_t^N),\mu_t)]\\
&\le    \ff c{\ss N}  \big(1+\|\mu_t\|_2\big)+ \int_0^t \e^{\int_s^t K_1(r)\d r} K_2(s) \E[\W_1(\mu^N(X_s^N),\mu_s)] \d s,\ \ t\ge 0.\end{align*}
Combining this with $\|\mu_t\|_2\le h_t(\mu_0)$ due to \eqref{ES20} for $p=2$, we derive   \eqref{AB1}  from Gronwall's inequality, while \eqref{AB2} follows from \eqref{AB1} and \eqref{K0}.

\end{proof}

 \beg{lem}\label{LN1} Assume $(H_1)$-$(H_2)$. Then for any $N\in\mathbb N$, $\nu,\tt\nu \in \scr P_1(\Z_{+}^N)$ and $t>0$, 
\beq\label{Y1}  \W_1(P_t^{N*} \nu,  P_t^{N*} \tt\nu)  \le \e^{\int_0^t  (K_1+K_2)(s) \d s} \W_1(\nu,\tt \nu).\end{equation}
 
\end{lem} 

\beg{proof}  For the initial value $(X_0^N,\tt X_0^N)$ with 
$$\L_{X_0^N}=\nu,\ \ \L_{\tt X_0^N} =\tt\nu,\ \ \E[\rr_N(X_0^N,\tt X_0^N)]= \W_1(\nu,\tt\nu),$$
by Theorem \ref{T0} let 
$(X_t^N, \tt X_t^N)$ be the inhomogenous Markov process on $\Z_+^N\times \Z_+^N$ generated by $\{\tt L_t^{N}: t\ge 0\}$, where
\beg{align*}   \tt L_t^{N} f(x,y) :=  \sum_{l=1}^N   &\Big[ \big(b_t(x_l,\mu^N(x))\land b_t(y_l,\mu^N(y))\big) \big(f(x+e_l,y+e_l)-f(x,y)\big) \\
&+ \big(a_t(x_l,\mu^N(x))  \land a_t(y_l,\mu^N(y)) \big) \big(f(x-e_l,y-e_l)-f(x,y)\big)\\ 
&+  \big(b_t(x_l,\mu^N(x))- b_t(y_l,\mu^N(y))\big)^+ \big(f(x+e_l,y)-f(x,y)\big)  \\
&+ \big(b_t(x_l,\mu^N(x))- b_t(y_l,\mu^N(y))\big)^- \big(f(x,y+e_l)-f(x,y)\big) \\
& +\big(a_t(x_l,\mu^N(x))- a_t(y_l,\mu^N(y))\big)^+ \big(f(x-e_l,y)-f(x,y)\big) \\
&+ \big(a_t(x_l,\mu^N(x))- a_t(y_l,\mu^N(y))\big)^- \big(f(x,y-e_l)-f(x,y)\big)\Big].\end{align*}
  Then as in \eqref{CP}, we have
\beg{align*}&\tt L_t^{N} \rr_N (x,y) \le \sum_{l=1}^N\big[K_1(t)|x_l-y_l|+ K_2(t)   \W_1(\mu^N(x),\mu^N(y))\big]\\
&= K_1(t)\rr_N(x,y)+ N K_2(t)   \W_1(\mu^N(x),\mu^N(y)),\ \ t\ge 0,\ x,y\in \Z_+^N.\end{align*}
So,
  \beg{align*}  
 &\W_1(P_t^{N*} \nu,  P_t^{N*} \tt\nu)   \le  \E[\rr_N(X_t^N,\tt X_t^N)]\\
 &\le \E[\rr_N(X_0^N,\tt X_0^N)]+\!\!\int_0^t\! K_1(s)\E[\rr_N(X_s^N,\tt X_s^N)]\d s+\!N\!\!\int_0^t\! K_2(s)\E[\W_1(\mu^N(X_s^N),\mu^N(\tt X_s^N)) ]\d s\\
&\le \E[\rr_N(X_0^N,\tt X_0^N)]+\int_0^t K_1(s)\E[\rr_N(X_s^N,\tt X_s^N)]\d s+\int_0^t K_2(s)\E[\rr_N(X_s^N,\tt X_s^N)]\d s\\
& \le \E[\rr_N(X_0^N,\tt X_0^N)]+\int_0^t (K_1(s)+K_2(s))\E[\rr_N(X_s^N,\tt X_s^N)]\d s,\ \ t\ge 0.
\end{align*} 
By Gronwall's inequality, we derive \eqref{Y1}.  
 \end{proof} 
 \section{Appendix:   inhomogenous  jump processes}\label{app}
 
%
 
 Let $E$ be a locally compact Polish space with Borel $\si$-algebra $\EE$, and let 
 $$q: [0,\infty)\times E\times \EE\to [0,\infty)$$ such that
 for each $t\ge 0$,   $q_t$ gives a stable conservative $q$-pair, i.e. 
 for every $x\in   E,$ $ q_t(x,\cdot)$ is a finite measure on $E$ with $q_t(x,\{x\})=0;$
 and for each  $A\in \EE,$  $q_t(\cdot,A)$ is a measurable function on $E$. Define $q_t(x):=q_t(x,E\backslash\{x\})$. Let $\B_b(E)$ be the class of bounded measurable  functions on $E$.

 \beg{enumerate} \item[$(A_1)$] For each $t\ge 0$, there exists a unique homogenous jump process with transition probability kernel $\{P^{(t)}_s(x,\d y)\}_{s\ge 0, x\in E}$  satisfying 
 the Kolmogorov equations 
 $$\pp_s P^{(t)}_sf= L_t P_s^{(t)}f= P_s^{(t)} L_t f,\ \ f\in \mathscr{B}_b(E),\ s\ge 0,$$ 
 where $P_s^{(t)}f(x):=\int_E f(y)P^{(t)}_s(x,\d y)$ and 
 $$L_tf(x):= \int_{E} \big(f(y)-f(x)\big) q_t(x,\d y),\ \ \ f\in C_0(E).$$
 \item[$(A_2)$] The $q$-pairs $(q_{t}(x,\cdot))_{t\ge 0,x\in E}$ satisfies 
 $$\lim_{\vv\downarrow 0} \sup_{s,t\le T, |t-s|\le \vv}   \|q_t(x,\cdot)-q_s(x,\cdot) \|_{var}  =0,\ \ T>0, $$locally uniform in $x\in E$,
 where $\|\cdot\|_{var}$ is the total variation distance. 
 \item[$(A_3)$] There exist a measurable function  $V\ge 1$ on $E$ with set $\{x:V(x)\le c\}$ being compact for any $c>0$, $K\in C([0,\infty); (0,\infty))$ and a constant $\theta>1$ such that 
 $$q_t(x)\le K(t)V(x),\ \   L_t V^\theta (x)\le K(t) V^\theta (x),\ \ \ \ t\ge 0, x\in E.$$
 \end{enumerate}

 \beg{thm}\label{T0} Assume $(A_1)$-$(A_3)$. There exists a unique inhomogeneous Markov process generated by $(L_t)_{t\ge 0}$, which is a family of transition  probability kernels 
 $\{P_{s,t}(x,\cdot):\ 0\le s\le t<\infty, x\in E\}$ on $E$, such that
 $$P_{s,t}f(x):= \int_E f(y)P_{s,t}(x,\d y),\ \ \ t\ge s\ge 0,\ f\in \B_b(E)$$
 satisfies    the semigroup property $P_{s,t}=P_{s,r}P_{r,t} $ for $0\le s\le r\le t$, and the  
 Kolmogorov equations 
 \beq\label{K1} \pp_sP_{s,t} f:= \lim_{\vv\downarrow 0} \ff{P_{s+\vv,t}f-P_{s,t}f}\vv= - L_s P_{s,t}f,\ \ \ t>s\ge 0,\ f\in \B_b(E),\end{equation}
 \beq\label{K2} 
  \pp_t P_{s,t}f:= \lim_{\vv\downarrow 0} \ff{P_{s,t+\vv}f-P_{s,t}f}\vv  =  P_{s,t}L_t f,\ \ \  t>s\ge 0,\ f\in\B_b(E).\end{equation} 
  
\end{thm}

\beg{proof} 
(a) By $(A_1)$,   for any $n\in \mathbb N$ and $k\in \Z_+$, let $\{P^{n,k}_t(x,\cdot)\}_{t\ge 0, x\in E}$ be the  unique Markov transition kernel   generated by
$L_{k 2^{-n}}$ satisfying the Kolmogorov equations.  For any $t\ge 0$, let
$$k_n(t):= \sup\{k\in \Z_+:\ k 2^{-n}\le t\}.$$
We now construct the transition probability kernel $\{P_{s,t}^{(n)}(x,\d y)\}$    
  for any $0\le s\le t<\infty, x\in E$ and $A\in \EE$   piecewisely in   time intervals $\{[k2^{-n}, (k+1)2^{-n}): k\in \Z_+\}$:   
\beg{align*}  &P_{s,t}^{(n)}(x,A):=   P_{t-s}^{n,k_n(s)}(x,A) \   \text{if}\  k_n(t)-k_n(s)=0; \\
&\!\!:= \int_{E} P_{t-k_n(t)2^{-n}}^{n,k_n(t)} (y_1,A)  P_{(k_n(s)+1)2^{-n}-s}^{n,k_n(s)} (x,\d y_1) \   \text{if}\ k_n(t)-k_n(s)=1; \\
&\!\!:= \!\! \int_{E^{l} }\!\! P_{t-k_n(t) 2^{-n}}^{n,k_n(t) } (y_{l},A)  P_{(k_n(s)+1)2^{-n}-s}^{n,k_n(s)} (x,\d y_1) 
\!\prod_{i=1}^{l-1} P_{2^{-n}}^{n,k_n(s)+i}(y_i,\d y_{i+1}) \   \text{if}\ k_n(t)\!-\!k_n(s)=l\ge 2.\end{align*}   
By $(A_1)$, it is  easy to check that $\{P^{(n)}_{s,t}(x,\cdot)\}_{0\le s\le t<\infty, x\in E}$ is an inhomogeneous Markov process generated by  
$$L_t^{(n)}  := L_{k_n(t)2^{-n}},\ \ \ t\ge 0,$$  in the sense that 
$$P_{s,t}^{(n)} f(x):= \int_E f(y)P_{s,t}^{(n)}(x,\d y),\ \ \ t\ge s\ge 0,\ f\in \B_b(E)$$
 satisfies    the semigroup property $P_{s,t}^{(n)}=P_{s,r}^{(n)}P_{r,t}^{(n)} $ for $0\le s\le r\le t$, and the  
 Kolmogorov equations 
 \beq\label{K1n} \pp_sP_{s,t}^{(n)} f = - L_s^{(n)} P_{s,t}^{(n)}f,\ \ \ t>s\ge 0,\ f\in \B_b(E),\end{equation}
 \beq\label{K2n}  \pp_t P_{s,t}^{(n)}f  =  P_{s,t}^{(n)}L_t^{(n)} f,\ \ \  t>s\ge 0,\ f\in\B_b(E).\end{equation} 
Below we show that when $n\to\infty$,  $P_{s,t}(x,\d y):= \lim_{n\to\infty} P_{s,t}^{(n)}(x,\d y)$ exists and is the unique jump process generated by $(L_t)_{t\ge 0}$ as claimed.

(b)    For any $f\in \B_b(E)$, let
$$P_{s,t}^{(n)} f(x):= \int_{E} f(y) P_{s,t}^{(n)}(x, \d y),\ \ 0\le s\le t<\infty, x\in E.$$ 
For any $x\in E$ and fixed $T>0$, define $$\varepsilon_m(x):=\sup\limits_{n>m, t\in [0,T+1]}(|b_{k_n(t)2^{-n}}(x)-b_{k_n(t)2^{-m}}(x)|+|a_{k_n(t)2^{-n}}(x)-a_{k_n(t)2^{-m}}(x)|),\ m\ge 1.$$
By $(A_2)$ we know that $0\le \vv_n(x)\to 0$ as $n\to \infty$. By $(A_3)$,  for any $x\in E$, there exists a constant $c>0$ such that   
\beq\label{*Y} \sup_{t\in [0,T+1]}  |(L_t^{(n)}-L_t^{(m)})f(x)| \le 2\|f\|_\infty \vv_{m}(x)\le c \|f\|_\infty V(x),\ \ n\ge m\ge 1.\end{equation} 
Next, by \eqref{K1n} and \eqref{K2n} we obtain the Duhamel formula 
$$(P_{s,t}^{(n)}-P_{s,t}^{(m)})f=\int_{s}^{t}P_{s,u}^{(n)}(L_{u}^{(n)}-L_u^{(m)})P_{u,t}^{(m)}f\d u,\ \ f\in \B_b(E).
$$
So, 
\begin{equation}\label{du}\sup_{0\le s\le t\le T}|(P_{s,t}^{(n)}-P_{s,t}^{(m)})f(x)|\le  2T\|f\|_{\infty}\sup_{n\ge 1, 0\le s\le t\le T}  P_{s,t}^{(n)} \vv_{m}(x).\end{equation}
By $(A_3)$ we have 
\beq\label{0} R:=\sup_{n\ge 1, 0\le s\le t\le T} P_{s,t}^{(n)} V^{\theta}(x) <\infty,\end{equation} 
where $\theta >1$. For any $\varepsilon>0$, let $M$ be the constant such that $cR/M^{\theta-1}<\varepsilon/2$. Since by $(A_2)$,  $\vv_n(x)\to 0$ uniformly in $x\in\{V\le M\}$ as $n\to\infty$, we know that there exists $m_0=m_0(\varepsilon)>0$ such that  $\vv_{m}(x)<\epsilon/2$ for $x\in \{V\le M\}$. Then
\begin{align*}
\sup_{n\ge 1, 0\le s\le t\le T}P_{s,t}^{(n)} \vv_{m}(x)&\le\sup_{n\ge 1, 0\le s\le t\le T}P_{s,t}^{(n)} (\vv_{m}(x) 1_{\{V\ge M\}})+\sup_{n\ge 1, 0\le s\le t\le T}P_{s,t}^{(n)} (\vv_{m}(x) 1_{\{V\le M\}})\\&\le cR/M^{\theta-1}+\sup_{n\ge 1, 0\le s\le t\le T}P_{s,t}^{(n)} (\vv_{m}(x) 1_{\{V\le M\}})\\&\le \varepsilon/2+\varepsilon/2=\varepsilon.
\end{align*}
Combining this  with \eqref{du}, we derive
\begin{equation}\label{1}\lim_{n\to\infty}\sup_{0\le s\le t\le T} \|P_{s,t}^{(n)}(x,\cdot)- P_{s,t} (x,\cdot) \|_{var} =0,\ \ x\in E, T>0 \end{equation}
  for  some Markov transition probability kernel $\{P_{s,t}(x,\cdot)\}_{0\le s\le t<\infty, x\in E}$, which is   an inhomogeneous  Markov jump process generated by $(L_t)_{t\ge 0}$ as wanted.
  Indeed, for any $f\in\B_b(E), x\in E$ and $0\le s<t\le T,$ we have
  $$\ff{P_{s+\vv,t}^{(n)}f(x)- P_{s,t}^{(n)} f(x)}{\vv} =-\ff 1 \vv\int_{s}^{s+\vv} L_r^{(n)} P_{r,t}^{(n)}f(x)\d r,\ \ \vv\in (0, t-s).$$
  By $(A_1)$ and \eqref{1}, when $n\to\infty$ this gives
$$\ff{P_{s+\vv,t}f(x)- P_{s,t}f(x)}{\vv} =-\ff 1 \vv\int_{s}^{s+\vv} L_rP_{r,t}f(x)\d r,\ \ \vv\in (0, t-s),$$
so that \eqref{K1} holds. Similarly, by
\beq\label{o} \ff{P_{s,t+\vv}^{(n)}f(x)- P_{s,t}^{(n)} f(x)}{\vv} =\ff 1 \vv\int_{t}^{t+\vv}   P_{s,r}^{(n)} L_r^{(n)} f(x)\d r,\ \ \vv\in (0, T-t).\end{equation} 
 By $(A_1)$ and $(A_2)$ we have $L_r^{(n)} f \to L_rf$ as $n\to\infty$, and there exists a constant $c>0$ such that 
 \beq\label{r}\sup_{r\in [0,T+1],n\ge 1} |L_r^{(n)} f|\le c V.\end{equation}
 Combining this with \eqref{0} and \eqref{1}, we may letting $n\to\infty$ in \eqref{o} to derive 
 $$ \ff{P_{s,t+\vv} f(x)- P_{s,t} f(x)}{\vv} =\ff 1 \vv\int_{t}^{t+\vv}   P_{s,r}  L_r  f(x)\d r,\ \ \vv\in (0, T-t).$$
 So, \eqref{K2} holds. Indeed, via \eqref{0} it suffices to show
 \beq\label{LB}  \lim_{n\to\infty} |P_{s,r}^{(n)} L_r^{(n)} f(x)-P_{s,r} L_r f(x)|=0.\end{equation}  To this end,  we note that 
 \begin{align*}
 	&|P_{s,r}^{(n)} L_r^{(n)} f(x)-P_{s,r} L_r f(x)|\\
	&\le |P_{s,r}^{(n)} L_r^{(n)} f(x)-P_{s,r} L_r^{(n)} f(x)|+	|P_{s,r} L_r^{(n)} f(x)-P_{s,r} L_r f(x)|\\&=:I_1(n)+I_2(n).
 \end{align*}
For any $N\in\mathbb{Z}_{+}$, we have 
\begin{align*}
	I_1(n)&\le  P_{s,r}^{(n)}(|L_r^{(n)}f| 1_{\{V\ge N\}})(x)+P_{s,r}^{(n)} (|L_r^{(n)} f| 1_{\{V\ge N\}})(x)\\&+|P_{s,r}^{(n)} (L_r^{(n)} f1_{\{V\le N\}})(x)-P_{s,r} (L_r^{(n)} f1_{\{V\le N\}})(x)|\\
	&\le   2cP_{s,r}^{(n)} (V 1_{\{V\ge N\}})(x)+|P_{s,r}^{(n)} (1_{\{V\le N\}}L_r^{(n)} f )(x)-P_{s,r} (1_{\{V\le N\}} L_r^{(n)}f )(x)|\\
	& \le    \frac{2c(P_{s,t}^{(n)}V^{\theta})(x)}{N^{\theta-1}}+|P_{s,r}^{(n)} (1_{\{V\le N\}}L_r^{(n)} f)(x)-P_{s,r} (1_{\{V\le N\}} L_r^{(n)} f )(x)|.
\end{align*}
Since  $\theta>1$, by \eqref{0}, \eqref{1}, \eqref{r},  we may first letting $n\to\infty$ then $N\to\infty$ to derive 
 \begin{align*}
	 \lim_{n\to\infty} I_1(n) =0.
\end{align*} 
 Moreover, by  \eqref{r} and dominated convergence theorem we have $I_2(n)\to 0$ as $n\to\infty$. Therefore,  \eqref{LB} is verified. 

(c) If there are two inhomogenous Markov processes generated by $(L_t)_{t\ge 0}$ with semigroups $(P_{s,t}^i)_{i=1,2}$ satisfying the Kolmogorov equations,  then  
follows from  the Duhamel formula
$$P_{s,t}^1f = P_{s,t}^2 f+\int_s^t P_{s,r}^1 (L_r-L_r) P_{r,t}^2 f\d r=0,\ \ \ f\in \B_b(E),$$ 
which is  implied by  the Kolmogorov equations. So, the uniqueness is confirmed. 
 \end{proof}
\section*{Declarations}
\noindent\textbf{Conflict of interest}\ {The authors declare that they have no conflict of interest that might influence the results
	reported in this paper.}
	
\noindent\textbf{Data Availability}\ {Data sharing is not applicable to this article as no new data were created or analyzed in this study.}

\noindent\textbf{Author Contributions} \ {All authors contributed to this work.}

\noindent\textbf{Acknowledgement} \ {The authors would like to thank the referee for corrections and very helpful suggestions.}


\end{document}